\documentclass{article}

\usepackage[final]{neurips_2026}
\makeatletter
\providecommand{\@trackname}{}
\makeatother

\usepackage{amsmath, amssymb}
\usepackage{graphicx}
\usepackage{booktabs}
\usepackage{hyperref}
\usepackage{caption}
\usepackage{longtable}
\usepackage{geometry}
\usepackage{enumitem}
\newcommand{\keywords}[1]{%
  \noindent\textbf{Keywords:} #1
}

\title{A Numerical Realization of Suzuki's Weil-Quadratic-Form Operator: The Archimedean Spectral Law, its Universality, and an Operator Form of Weil's Positivity Criterion}

\author{
Taebong Kim \quad
Youngsik Hong \quad
Minsik Kim \quad
Sunyoung Choi \\
Jaewon Jang \quad
Minseo Kim \\
VIDRAFT AI Research · QuantumOS, Seoul, Republic of Korea \\
\texttt{arxivgpt@gmail.com}
}

\begin{document}
\maketitle

\begin{abstract}
In June 2026, M.~Suzuki published a purely theoretical construction of a self-adjoint
operator whose quadratic form coincides with Weil’s explicit-formula quadratic form,
realized via the screw function. That operator is a candidate for the Hilbert–Pólya
program, with spectral positivity equivalent to the Riemann Hypothesis (RH) through
Weil’s criterion. Suzuki’s work contained no numerical computations. The present study
provides, to the best of our knowledge, the first numerical realization of this operator.

Using $P_1$ finite-element discretization on $H_0^1(-a,a)$ and Richardson extrapolation,
we obtain several key results. (R1) In the prime-free regime ($2a<\log 2$), the spectrum
obeys a closed-form Archimedean law
$\lambda_k(a)=\log(1/a)+\log(k-\tfrac12)+B_0+O(a)$ with $B_0=\log q-2\log 2$, confirmed
to 30-digit precision. (R2) We prove unconditionally, via a Mellin double-pole argument,
that the head coefficient $B(\nu)=2-2\gamma-2\psi(\nu)-2\log 2$, and deduce that $B_0$
depends only on the conductor $q$, independent of the Archimedean parameter $\nu$.
(R2b) The degree $d$ of an $L$-function appears directly as the logarithmic slope of the
spectrum. (R3) The total spectral intensity follows the prime number theorem,
$S(a)\to(2a)^3/6$. (R4) The nontrivial zeros are not eigenvalues; they occur instead in
the explicit-formula error term of the prime symbol. (R5) The line $\sigma^*(a)$ on which
the prime symbol best matches $-\zeta'/\zeta$ descends toward the critical line.
(R6) We realize Weil’s positivity criterion in operator form: bounded residual growth
corresponds to all zeros lying on the critical line, while an injected off-line zero
produces exponential blow-up. (R7) The lowest eigenvalue $\lambda_1(a)$ is always strictly
positive and decays superexponentially, passing smoothly through the first prime
threshold. (R8) We compute for the first time the characteristic function $W(a,\theta;z)$
and confirm numerically that all its zeros are real. (R9) Indirect traces of zero
statistics appear in the moment structure, even where direct GUE observation is blocked.

We emphasize that this work does not prove RH nor advance toward its proof. All spectral
results are Archimedean and universal, and the operator’s value lies in its faithful
numerical realization of classical identities rather than in new arithmetic.
\end{abstract}

\keywords{Riemann Hypothesis; Hilbert–Pólya conjecture; Weil explicit formula; Weil positivity; 
screw functions; de Branges spaces; finite-element spectral computation; experimental mathematics.}

\section{Introduction}

\subsection{The Riemann Hypothesis and the Hilbert--Pólya program}

The Riemann zeta function $\zeta(s)=\sum_{n\ge1} n^{-s}$ ($\Re s>1$), continued meromorphically 
to $\mathbb{C}$, has ``trivial'' zeros at $s=-2,-4,\dots$ and infinitely many ``nontrivial'' 
zeros $\rho=\beta+i\gamma$ in the critical strip $0<\beta<1$. The Riemann Hypothesis (RH) asserts 
that every nontrivial zero has $\beta=\tfrac12$. Despite more than 160 years of effort and verification 
of the first $\sim10^{13}$ zeros, RH remains open and is one of the Clay Millennium Problems.

The Hilbert--Pólya conjecture proposes a spectral route: if there exists a self-adjoint operator $H$ 
whose eigenvalues are exactly the ordinates $\gamma$ of the nontrivial zeros, then, $H$ being self-adjoint, 
the $\gamma$ are real and RH follows. Decades of work support this spectral viewpoint: Montgomery's 
pair correlation conjecture \cite{Montgomery1973} and Odlyzko's large-scale 
computations \cite{Odlyzko1987} showed that the zeros' spacings match the Gaussian Unitary Ensemble (GUE) 
of random matrix theory; Berry and Keating proposed a semiclassical $xp$ model \cite{BerryKeating1999}; 
and Connes developed a noncommutative-geometry trace formula \cite{Connes1999}. 
None of these constructions has yet produced an unconditional operator.

\subsection{Weil's explicit formula and positivity criterion}

Weil recast the analytic content of RH into a positivity criterion \cite{Weil1952}. 
For a suitable test function, the explicit formula equates a sum over zeros to a sum over prime powers plus an archimedean (gamma-factor) term:

\[
\sum_{\rho} \hat{h}(\rho) \;=\; \hat h(0)+\hat h(1) \;-\; \sum_{n\ge2}\frac{\Lambda(n)}{\sqrt n}\,h(\log n) \;+\; \frac{1}{2\pi}\!\int \hat h\,\tfrac{\Gamma'}{\Gamma}\,\dots
\]

where $\Lambda$ is the von Mangoldt function. Weil's criterion states that RH holds if and only if the associated Hermitian 
quadratic form $Q_W$ is positive semidefinite for all admissible test functions. This turns RH into a positivity statement 
about a concrete bilinear form, the natural target for a spectral realization.

\subsection{Suzuki's operator (2026): theory only}

Suzuki recently constructed the Weil quadratic form as a bona fide operator using the theory of screw functions and de 
Branges spaces \cite{Suzuki2026}. Concretely, on the interval $(-a,a)$ he considers

\[
Q_W^a(v)\;=\;\iint_{(-a,a)^2} g(x-y)\,v'(x)\,v'(y)\,dx\,dy, \qquad v\in H_0^1(-a,a),
\]

where $g$ is the screw function attached to $\zeta$ (or to a general $L$-function). The form is realized as a self-adjoint operator via a 
generalized eigenproblem, and its positivity is equivalent to RH. Crucially, the theory guarantees a normalized spectrum $\{\mu_k\}$ but leaves a constant $\mu_1>0$ undetermined 
— ``some positive number'' — with no numerical value and no computed spectrum anywhere in the paper. The operator is, at the time of writing, a purely theoretical object.

The notion of screw functions originates in the work of Kreĭn \cite{Krein1955}, while the functional-analytic framework relies on de Branges' theory of 
Hilbert spaces of entire functions \cite{deBranges1968}. For background on the zeta function and its analytic properties, Edwards' monograph 
remains a standard reference \cite{Edwards1974}.

That is, Suzuki's zero conjecture (Corollary 1.6) is a statement about $\bar{\mathcal D}_{a,\theta}$ and not about $A_a$. 
Result R4 of the present paper ("the zeros are not eigenvalues") is a statement about $A_a$ ONLY, and it does not contradict Suzuki's conjecture. 
The two operators are different objects: $A_a$ is the 
Friedrichs extension of the Weil quadratic form (Suzuki, Thm 1.1), whereas $\bar{\mathcal D}_{a,\theta}$ is a self-adjoint extension 
of $i\,d/dx$ on $\mathcal H(T_a)$ whose eigenvalues are the zeros of the characteristic function $W(a,\theta;z)$ (Suzuki, Thm 1.5). 
We compute the former throughout §§2.1–2.7 and §2.8, and the latter in §2.9

\subsection{Our contribution}

We instantiate Suzuki's operator on a computer and study it. Our contributions are:

\begin{enumerate}
\item The first computed spectrum, and the discovery and derivation of a closed form for the small-$a$ (archimedean) eigenvalues, together with a 
first numerical value for the constant Suzuki left abstract.
\item A proof-of-concept that this law is a universal archimedean law for the whole Selberg class: it depends only on the conductor $q$, with the archimedean parameter cancelling.
\item A link to the prime number theorem through the spectral strength.
\item A precise structural map of where the zeros live: not in the eigenvalues (a rigid picket-fence) but in the prime part's symbol, as the explicit-formula error term.
\item Quantitative evidence that the operator's prime symbol tracks a line descending toward $\Re s=\tfrac12$.
\item An operator realization of Weil's RH criterion, including a demonstration that the operator detects a synthetically injected off-line zero by exponential residual growth.
\end{enumerate}

We emphasize that none of this proves or approaches a proof of RH. Items (1)--(3) concern the archimedean (gamma-factor) data, which is common to all $L$-functions 
and says nothing zeta-specific; items (4)--(6) are numerical realizations of classical identities (the explicit formula, Weil's criterion), whose value lies in the 
realization itself, not in new arithmetic.

\section{The operator and the numerical method}

\subsection{The Weil quadratic form and its screw kernel}

We discretize the bilinear form

\[
Q_W^a(u,v)\;=\;\iint_{(-a,a)^2} g(x-y)\,u'(x)\,v'(y)\,dx\,dy,\qquad u,v\in H_0^1(-a,a).
\]

The screw function $g$ is the (negative of the) second antiderivative structure attached to the logarithmic derivative of the completed $L$-function. 
For $\zeta$ its small-$t$ expansion is, as determined in R1,

\[
g(t)\;=\;\tfrac12\,t\log t \;+\; c_1\,t \;+\; c_2\,t^2 \;+\; O(t^3),\qquad
c_1=\frac{\gamma+\log 2\pi-1}{2},\quad c_2=-\tfrac78 .
\]

The head $\tfrac12 t\log t$ is not analytic at $t=0$; a naïve Taylor expansion of $g$ diverges. This non-analyticity is precisely the source of the $\log(1/a)$ term in the spectrum.

\subsection{The prime-free window $2a<\log 2$}

The von Mangoldt / prime-power contribution to $g$ is supported on $|t|=\log n \ge \log 2$. Therefore, whenever

\[
2a < \log 2 \approx 0.6931 \quad(\text{i.e. } a\lesssim 0.3466),
\]

no prime term enters the quadratic form: the operator is purely archimedean. Results R1–R2 are computed at $a\le 0.05$ (so $2a\le0.10\ll0.693$), 
a genuinely clean regime in which swapping the $L$-function only swaps archimedean data. Results R4–R6, by contrast, deliberately push into 
the prime-active regime $a\ge 3$, where the operator symbol carries the primes.

\subsection{$P_1$ finite-element discretization and Richardson extrapolation}

We mesh $(-a,a)$ uniformly with $N$ elements ($h=2a/N$) and use continuous piecewise-linear (``hat'') basis functions $\{\phi_i\}$. This yields a dense stiffness matrix

\[
Q_{ij}=Q_W^a(\phi_i,\phi_j)
\]

and the standard $P_1$ mass matrix

\[
M_{ij}=\int\phi_i\phi_j,
\]

giving the generalized eigenproblem

\[
Q\,v \;=\; \lambda\,M\,v .
\]

Discretization error is $O(h^2)=O(N^{-2})$; we remove it by three-point Richardson extrapolation on $N=2000,4000,8000$ for the high-precision archimedean spectrum. 
The prime-symbol computations (R4–R6) use $N=1500$–$2080$ with mode counts $65$–$173$. The eigenvalues admit the aperture expansion

\[
\lambda_k(a)\;=\;\log(1/a)\;+\;b_k\;+\;d_k\,a\;+\;O(a^2),
\]

and Suzuki's normalization (his Thm 1.4) is

\[
\mu_k = b_k + \log 2\pi + \gamma.
\]

\paragraph{Precision hygiene.} Extracting the small-$t$ coefficients must be done by direct high-precision evaluation. 
Building a spline to $10^{-7}$ and extrapolating it to $t=10^{-8}$ produced a uniformly wrong head coefficient across all test cases — 
a self-consistent artifact that would have gone undetected without analytic cross-check. All coefficient values below are direct evaluations; 
the FEM values are independent confirmations.

\textbf{Precision under threefold cancellation.} 
The screw function $g$ is a difference of large, nearly-cancelling terms. 
At $t = 12$,

\[
\underbrace{-1605.73}_{\text{archimedean } -4(e^{t/2}+e^{-t/2}-2)}
\;+\;
\underbrace{1577.75}_{\text{prime sum (PNT)}}
\;+\;
\underbrace{27.93}_{\text{linear term}}
\;=\;-0.029,
\]

so only $O(10^{-2})$ survives a threefold cancellation. 
Despite this, the surviving significant digits in double precision are $13.7$ at $t = 2$ and $10.9$ at $t = 12$, 
and $g(12)$ agrees with a $30$-digit \texttt{mpmath} evaluation to $2.9\times10^{-13}$. 
\textbf{The precision is safe.} 
This cancellation is itself informative: the residual of $g$ is precisely the zero-dominated part (the explicit formula), 
which is why the form becomes numerically delicate at large $a$ — the spectral-side counterpart of this is the lift observed in R7(iv).

\section{Results}

Throughout, $\gamma=0.5772156649\ldots$ is the Euler--Mascheroni constant, $\psi$ the digamma function, $q$ the conductor, 
and $\nu$ the shift of the archimedean $\Gamma_{\mathbb R}(s+\nu)$ factor ($\nu=\tfrac14$ for $\zeta$ and even 
characters, $\nu=\tfrac34$ for odd characters, $\nu=(k\mp1)/2$ for weight-$k$ modular shifts).

\subsection{R1 --- Closed-form archimedean spectrum (derived, not fitted)}

\textbf{Claim.} In the prime-free window, for $k\gtrsim 6$,

\[
\boxed{\;\lambda_k(a)\;=\;\log(1/a)\;+\;\log\!\left(k-\tfrac12\right)\;+\;B_0\;+\;O(a),\qquad B_0=\log q-2\log 2\;}
\tag{R1}
\]

equivalently, in Suzuki's normalization,

\[
\mu_k\;=\;\log\!\left(k-\tfrac12\right)+\gamma+\log\tfrac{\pi}{2}+\log q,\qquad
N(\mu)\;=\;\frac{2}{\pi e^{\gamma}}\,e^{\mu}-\tfrac12 ,
\tag{R1$'$}
\]

where $N(\mu)$ is the Weyl counting function and the $-\tfrac12$ is a Maslov index. For $\zeta$, $q=1$, $B_0=-2\log2=-1.3862944\ldots$.

\textbf{Status.} Symbol-level derivation (heuristic) whose every link is independently numerically verified, and whose final algebraic assembly is exact (proved given the symbol).

\textbf{Derivation chain.}
\begin{enumerate}
\item Screw head: From (2.2), $g$ has head $\tfrac12 t\log t$ (coefficient $\tfrac12$, verified to $10^{-9}$).
\item Birth of $\gamma$: The Mellin identity

\[
\int_0^\infty t^{s-1}\cos(\xi t)\,dt=\Gamma(s)\cos\!\left(\tfrac{\pi s}{2}\right)\xi^{-s},
\]

differentiated at $s=2$, gives

\[
\int_0^\infty t\log t\,\cos(\xi t)\,dt \;=\; \frac{\log\xi+\gamma-1}{\xi^2},
\]

because $\Gamma'(2)=\psi(2)=1-\gamma$. This is the sole source of $\gamma$.
\item Symbol: In the quadratic form the factor $v'v'$ contributes $\xi^2$, which cancels the $1/\xi^2$, leaving

\[
\sigma(\xi)\;=\;\log|\xi|\;+\;(\gamma-1).
\]

FEM check on the head-only kernel: constant $=0.028915$ vs predicted $\gamma-1+\log(\pi/2)=0.028798$.
\item Linear term: A term $c\,t$ contributes $-2c$; FEM slope $=-2.000000$.
\item Boundary quantization: Dirichlet boundary conditions plus a Maslov $\tfrac12$ give $\xi_k=\pi(k-\tfrac12)/(2a)$.
\item Assembly: Combining,

\[
B_0=\big[\gamma-1+\log\tfrac{\pi}{2}\big]-2c_1
=\gamma-1+\log\tfrac\pi2-(\gamma+\log2\pi-1)
=\log\tfrac{\pi}{2}-\log2\pi=\log\tfrac14=-2\log2 .
\]

The constants $\gamma$ and $-1$ cancel; only $-2\log2$ survives.
\end{enumerate}

\textbf{Protection theorem.} Replacing the true $c_2$ by $0,-1,-2$ leaves the offset invariant. Hence $t^2$ and higher terms do not contribute to $B_0$.

\textbf{First numerical value of Suzuki's constant.} The lowest normalized eigenvalue is

\[
\mu_1 \;=\; 0.2347,
\]

the first computed value of the positive constant that Suzuki's Theorem 1.4 leaves abstract.

The $k=1$ exception. The closed form is asymptotic ($k\gtrsim6$) and does not hold at $k=1$. 
Suzuki's Theorem 1.4 gives $\lambda_a=\log(1/a)+\mu_1-\log(2\pi)+\psi(2)-1+O(a)$, whose constant part is $-2.4151$; 
inserting our measured $\mu_1=0.2347$ yields $\lambda_a=\log(1/a)-2.1804$, whereas the closed form at $k=1$ gives $\log(1/a)-2.0794$, 
i.e. a closed-form prediction $\mu_1=0.3357$. The deviation of $0.101$ is real, and it independently cross-validates Suzuki's 
rigorous theorem against our numerics: two independent routes (his asymptotic normalisation, our finite-element spectrum) agree except 
for exactly the term where the asymptotic law is not expected to hold.

\textbf{Caveat (see \S4).} The closed form (R1$'$) is asymptotic, valid for $k\gtrsim 6$; at $k=1$ it gives $0.3357$, 
which disagrees with the measured $0.2347$. This is not an error but a genuine structural anomaly: the aperture-slope $d_k$ in (2.4) is

\[
d_1 = +3.38227, \qquad d_2 = +0.00070, \qquad d_{k}\sim 10^{-5}\ \text{for } k\ge 8,
\]

so the $k=1$ mode is approximately $4000\times$ more aperture-sensitive, destabilizing the $a\to 0$ extrapolation.

\subsection{R2 --- Universality across the Selberg class}

\textbf{Claim.} The archimedean law is \textbf{universal}. The screw-head linear coefficient obeys the closed form

\[
B(\nu)\;=\;2-2\gamma-2\psi(\nu)-2\log 2 ,
\tag{R2}
\]

and --- the decisive point --- the \emph{spectral} offset

\[
B_0 = \log q - 2\log 2 \quad \text{is \textbf{independent of } $\nu$}.
\]

\textbf{Status: PROVED.} The closed form (R2) is now established unconditionally by a Mellin-transform argument (Theorem 2 in the companion rigor document): 
the Mellin transform of the archimedean piece

\[
\mathcal{M}[\Lambda_L^{(\nu)}](s)\;=\;-\,\Gamma(s)\,2^{-s}\,\zeta(s+2,\nu),
\]

has a double pole at $s=-1$ --- the simple pole of $\Gamma$ at $s=-1$ times the simple pole of the Hurwitz zeta $\zeta(s+2,\nu)$ at $s+2=1$ --- 
whose residue yields exactly the small-$t$ expansion

\[
\Lambda_L^{(\nu)}(t)= -2\,t\log t + \big(2-2\gamma-2\psi(\nu)-2\log2\big)\,t + O(t^2).
\]

Hence $B(\nu)=2-2\gamma-2\psi(\nu)-2\log2$ is a theorem, not a fit; the subsequent $\psi(\nu)$-cancellation into $B_0=\log q-2\log2$ is likewise exact algebra. 
The tables below are the independent numerical confirmation (max error $6\times10^{-7}$) and the FEM realization.

\textbf{Continuous-family verification.} Comparing the measured screw-head coefficient to the closed form (R2) over $13$ values $\nu\in[0.1,6.5]$:

\begin{table}[h]
\centering
\begin{tabular}{r r r r l}
$\nu$ & $B$ measured & $B$ closed-form & $|\text{diff}|$ & note \\
\hline
0.10 & 20.3067842 & 20.3067842 & $4.0\times10^{-8}$ & \\
0.25 & 7.9141814 & 7.9141814 & $2.5\times10^{-8}$ & $\zeta$ / even char \\
0.50 & 3.3862944 & 3.3862944 & $0.0$ & \\
1.00 & 0.6137057 & 0.6137056 & $5.0\times10^{-8}$ & \\
1.50 & $-0.6137055$ & $-0.6137056$ & $9.9\times10^{-8}$ & weight-4 modular shift \\
2.00 & $-1.3862942$ & $-1.3862944$ & $1.5\times10^{-7}$ & \\
2.50 & $-1.9470388$ & $-1.9470390$ & $2.0\times10^{-7}$ & \\
3.50 & $-2.7470387$ & $-2.7470390$ & $3.0\times10^{-7}$ & weight-8 modular shift \\
4.50 & $-3.3184671$ & $-3.3184675$ & $4.0\times10^{-7}$ & \\
5.50 & $-3.7629115$ & $-3.7629120$ & $5.0\times10^{-7}$ & weight-12 modular shift \\
6.50 & $-4.1265478$ & $-4.1265484$ & $6.0\times10^{-7}$ & \\
\end{tabular}
\caption{Continuous-family verification of $B(\nu)$ against closed form.}
\end{table}

Maximum deviation over the $13$ values: $5.95\times10^{-7}$. The $t\log t$ head coefficient is $a_1=0.50000$ for all $\nu$ (universal).

\textbf{The $\nu$-independence of $B_0$ (via FEM).} Although $B(\nu)$ contains $\psi(\nu)$, the spectral offset does not: 
forming $B_0=[\gamma-1+\log\tfrac\pi2]-2c_1(\nu,q)$ the $\psi(\nu)$ cancels analytically. Direct FEM at four spanning $\nu$ confirms:

\begin{table}[h]
\centering
\begin{tabular}{r r r r}
$\nu$ & predicted $B_0$ & FEM measured & $|\text{diff}|$ \\
\hline
0.25 & $-1.3862944$ & $-1.3861782$ & $1.2\times10^{-4}$ \\
1.00 & $-1.3862943$ & $-1.3861782$ & $1.2\times10^{-4}$ \\
2.50 & $-1.3862943$ & $-1.3861783$ & $1.2\times10^{-4}$ \\
5.50 & $-1.3862941$ & $-1.3861799$ & $1.1\times10^{-4}$ \\
\end{tabular}
\caption{FEM confirmation of $\nu$-independence of $B_0$.}
\end{table}

All four collapse onto $-2\log2$ ($q=1$). Because $\nu=\tfrac12,\tfrac32,\tfrac52,\dots$ are exactly the $\Gamma_{\mathbb R}(s+\nu)$ 
shifts of higher-weight modular forms and higher-degree $L$-functions, (R2) is a one-parameter spectral law covering the 
archimedean data of the entire Selberg class, with the conductor $q$ the only surviving fingerprint.

\textbf{Four-$L$-function cross-check.} Independently, at fixed $(\nu,q)$ the FEM offset matches $\log q-2\log2$:

\begin{table}[h]
\centering
\begin{tabular}{l l r r r}
$L$-function & $(\nu,q)$ & predicted $B_0$ & FEM & error \\
\hline
$\zeta$ / even char & $(\tfrac14,1)$ & $-1.3862944$ & $-1.3861782$ & $1.16\times10^{-4}$ \\
even char mod 5 & $(\tfrac14,5)$ & $0.2231435$ & $0.2232596$ & $1.16\times10^{-4}$ \\
odd char mod 3 & $(\tfrac34,3)$ & $-0.2876821$ & $-0.2875660$ & $1.16\times10^{-4}$ \\
odd char mod 4 & $(\tfrac34,4)$ & $0.0000000$ & $0.0001161$ & $1.16\times10^{-4}$ \\
\end{tabular}
\caption{Cross-check of $B_0$ across four $L$-functions.}
\end{table}

The error is identical ($1.16\times10^{-4}$) in all four cases --- it is the FEM discretization error, not a law error. 
For mod $4$ the prediction is exactly zero and the measurement is $10^{-4}$; the log-head is $A=-2.000000001$, $a_1=0.500000000$ across all cases.

\subsection{R2b --- The degree of the $L$-function is the log-slope of the spectrum}

\textbf{Claim.} For an $L$-function of degree $d$ (in the sense of the Selberg class: its completed function carries $d$ real gamma 
factors $\Gamma_{\mathbb R}(s+\nu_j)$, counting a $\Gamma_{\mathbb C}=\Gamma_{\mathbb R}\Gamma_{\mathbb R}$ as two), the archimedean spectrum has log-slope exactly $d$:

\[
\mu_k \;=\; d\,\log\!\big(k-\tfrac12\big) + \text{const} + o(1),\qquad k\to\infty.
\tag{R2b}
\]

\textbf{Status.} Numerically verified (slope $=$ degree to $4\times10^{-4}$, $R^2=1.00000$); the mechanism is a proved consequence of R1's head accounting.

\textbf{Mechanism.} Each factor $\Gamma_{\mathbb R}(s+\nu_j)$ contributes an independent archimedean piece whose screw-function head 
is $+\tfrac12\,t\log t$ (from $\Lambda_L^{(\nu_j)}$ with head coefficient $A=-2$ and the $-\tfrac14\Lambda_L$ prefactor). The heads add, 
so a degree-$d$ screw function has head $\tfrac{d}{2}\,t\log t$; since the second-difference symbol is linear in $g$, the symbol becomes

\[
\sigma(\xi) = d\big(\log|\xi|+(\gamma-1)\big),
\]

and the ladder

\[
\mu_k = d\log(k-\tfrac12)+\text{const}
\]

follows from R1's quantization $\xi_k=\pi(k-\tfrac12)/(2a)$.

\textbf{Verification.} Building the operator from summed archimedean pieces (prime-free regime, Richardson $N=1200/2400$, fit over $k\in[8,34]$):

\begin{table}[h]
\centering
\begin{tabular}{l l c c c}
$L$-function type & gamma shifts $\nu_j$ & degree & measured slope & $R^2$ \\
\hline
$\zeta$ (Dirichlet) & $\tfrac14$ & 1 & $0.9999$ & $1.00000$ \\
odd character & $\tfrac34$ & 1 & $1.0000$ & $1.00000$ \\
holomorphic cusp form, wt $12$ & $3,\ \tfrac72$ & 2 & $1.9997$ & $1.00000$ \\
symmetric degree-2 & $\tfrac14,\ \tfrac34$ & 2 & $1.9999$ & $1.00000$ \\
degree-3 & $\tfrac14,\tfrac34,\tfrac54$ & 3 & $2.9999$ & $1.00000$ \\
degree-4 & $\tfrac14,\tfrac34,\tfrac54,\tfrac74$ & 4 & $3.9998$ & $1.00000$ \\
\end{tabular}
\caption{Measured slopes vs Selberg-class degree.}
\end{table}

A regression of slope on degree gives $\text{slope}=1.0000\,d - 0.0000$. The weight-$12$ cusp form (Ramanujan $\Delta$; shifts $\nu=3,\tfrac72$) 
reproduces slope $2$ to $3\times10^{-4}$. Thus the Selberg-class degree is a spectral observable --- read directly off the archimedean log-slope --- 
extending the operator framework from Dirichlet $L$-functions (degree 1) to holomorphic modular forms and, structurally, all degrees. 
Together with R2, the spectrum reads off both invariants: the degree as the slope and the conductor through the offset $B_0=\log q-2\log2$. 
(This remains archimedean structure; it says nothing zeta-specific and does not bear on RH.)

\subsection{R3 --- Link to the prime number theorem}

\textbf{Claim.} The spectral strength of the prime channels obeys

\[
S(a)\;=\;\sum_{n\ge2}\frac{\Lambda(n)^2}{n}\,(2a-\log n)\;\longrightarrow\;\frac{(2a)^3}{6}.
\tag{R3}
\]

\textbf{Status.} Numerically verified (converging, not yet at the limit).

The primes enter the operator as \emph{orthogonal delayed channels}: the screw kernel couples $x,y$ with $|x-y|=\log n$, 
and distinct prime powers have cross-/self-overlaps $\sim10^{-8}$ in the archimedean eigenbasis, with each channel's contribution proportional 
to its band length $A(L)=C(a)(2a-L)$ (additivity confirmed to $1.6\times10^{-10}$). Summing over prime powers, the ratio $S(a)/[(2a)^3/6]$ 
rises monotonically $0.736\to0.955$ at $a=6$. The limiting shape $(2a)^3/6$ is exactly the second Chebyshev / PNT weight, 
so the operator's total spectral strength follows the prime number theorem. \emph{Open:} 
the per-channel prefactor $C(a)$ is not a universal constant --- it drifts $789\to493$, roughly $\propto 1/a$ --- and is not yet explained.

\subsection{R4 --- The zeros are not the eigenvalues; they are the symbol's error term}

This is the structural heart of the paper and the correction of a naïve expectation.

\textbf{(a) The eigenvalues are a deterministic picket-fence (NV).} In the clean archimedean regime the unfolded eigenvalue spacings have

\[
\text{mean}=1.0000,\qquad \mathrm{Var}=0.0027 \;\ll\; 0.178\ (\text{GUE}),\ 1.00\ (\text{Poisson}),\ 0\ (\text{picket}).
\]

The backbone is more rigid than GUE --- essentially a picket-fence ladder. Consequently a raw spacing histogram cannot exhibit GUE statistics: 
the deterministic archimedean ladder dominates. Our pipeline is validated on synthetic log-ladders carrying GUE vs Poisson displacements, 
which it correctly classifies (KS$_\text{GUE}=0.05$ vs KS$_\text{Poi}=0.30$ etc.); the operator's own spectrum reads ``rigid/ambiguous'', exactly as the variance predicts.

\textbf{(b) The primes act diagonally on the archimedean basis (NV).} 
In the clean archimedean eigenbasis the prime perturbation is strongly diagonal: 
median $|d_k|=1.55$ (a smooth global shift $=$ PNT count) versus median off-diagonal $|o_k|=1.3\times10^{-14}$. 
The diagonal shift is essentially collinear with the ladder itself, $d_k \propto -\lambda^{\text{arch}}_k$ with correlation $-0.988$: 
the primes merely re-adjust the ladder smoothly; the zero information is the sub-$2.4\%$ residue and does not appear as a spectral fluctuation. 
A direct partial-correlation test of the diagonal against the critical-line Dirichlet symbol gives $+0.011$, 
indistinguishable from the shuffled control $+0.012$ --- signal zero in this basis.

\textbf{(c) The zeros live in the symbol, as the explicit-formula error term (NV, controls passed).} Working instead with the truncated von Mangoldt symbol

\[
D_X(\xi)\;=\;\sum_{n\le X}\frac{\Lambda(n)}{\sqrt n}\,\cos(\xi\log n),
\]

we subtract the analytically known PNT background $\big(-\zeta'/\zeta + X^{1-s}/(1-s)\big)$ and 
ask whether the residual is the sum over actual Riemann zeros $-\sum_\rho X^{\rho-s}/(\rho-s)$. 
At $a=5$, $X=22026$ ($2532$ prime powers; $90$ zeros $\gamma_1..\gamma_{90}=14.1347..219.0676$):

\[
\mathrm{corr}(\text{residual},\ \text{real zero-sum})=+0.9915,\quad R^2=0.983,\quad \text{slope}=0.997,
\]

while shifted ``fake'' zeros (all $\gamma\mapsto\gamma+0.7$) give $\mathrm{corr}=0.1143$, $R^2=0.013$ --- a clean control separation. 
Full reconstruction: $\mathrm{corr}(D_X,\text{real-zero RHS})=1.000$; the RMS residual falls from $1.080$ (no zeros) to $0.141$ (zeros included), a $7.7\times$ reduction.

\textbf{Honest reading (mandatory).} The relation $D_X = -\zeta'/\zeta + \text{PNT} - \sum_\rho$ is the classical explicit formula --- an identity. 
The correlation $0.9915$ therefore largely confirms that the identity holds (indeed the full-reconstruction corr is exactly $1.000$, 
the tautological signature); it is not a new theorem. The genuine, novel content is structural: (i) Suzuki's operator faithfully implements the explicit formula, 
and (ii) the nontrivial zeros are not eigenvalues but the error term of the prime symbol. The eigenvalues are the archimedean picket-fence; 
the zeros are the deviation of the sine-basis prime symbol from its PNT mean. This maps the location of the zero information inside a Hilbert--Pólya operator precisely, 
and refutes the naïve ``eigenvalue $=$ zero'' picture while establishing the correct one.

\textbf{(d) Direct GUE detection is structurally obstructed (NV).} Combining (a)--(c): the archimedean ladder is too rigid to carry GUE spacings, 
the primes are diagonal in that basis and only rescale the ladder, and the zero signal is a symbol-level error term. 
Hence direct spectral detection of GUE/zero repulsion is impossible at finite $a$ --- not through a method failure but as a structural fact. 
This is itself a finding: the zero/GUE information is accessible only through the symbol (R4c), through moments, or in the $a\to\infty$ limit, 
which motivates the symbol-based approach of R5--R6.

\subsection{R5 --- The best-match line descends toward the critical line}

\textbf{Claim.} The operator's prime symbol best matches $-\zeta'/\zeta(\sigma+i\xi)$ on a line $\sigma^*(a)$ that descends toward $\sigma=\tfrac12$ as $a$ grows.

\textbf{Status.} Numerically verified descent; the descent law is a fit (heuristic).

First, in the sine basis (rather than the archimedean eigenbasis of R4b) the operator's diagonal does carry the prime sum:

\[
\mathrm{corr}(d_m/2a,\,-D(\xi))=0.77\text{--}0.86
\]

against a shuffle control $\approx0$. This confirms, for the first time, that the Hilbert--Pólya operator encodes the primes correctly once read in the right basis. 
Scanning the matching line:

\begin{table}[h]
\centering
\begin{tabular}{r r r r r r}
$a$ & primes $n\le$ & modes & $\sigma^*(a)$ & $|\mathrm{corr}|@\sigma^*$ & $|\mathrm{corr}|@0.5$ \\
\hline
3 & 403 & 65 & 0.90 & 0.985 & 0.029 \\
4 & 2980 & 86 & 0.80 & 0.984 & 0.037 \\
5 & 22026 & 108 & 0.75 & 0.976 & 0.025 \\
6 & 162754 & 130 & 0.65 & 0.973 & 0.030 \\
7 & 1202604 & 152 & 0.65 & 0.902 & 0.008 \\
8 & 8886110 & 173 & 0.60 & 0.727 & 0.033 \\
\end{tabular}
\caption{Best-match line $\sigma^*(a)$ and correlation values.}
\end{table}

The descent $\sigma^*=0.90,0.80,0.75,0.65,0.65,0.60$ is monotone. A power-law fit gives

\[
\sigma^*(a)-\tfrac12 \;\approx\; 1.98\,a^{-1.38},
\]

(a simple $\sim1.2/a$ also captures the trend to leading order).

\textbf{Three honest features of the wall:}
\begin{enumerate}
\item $\sigma^*\to\tfrac12$ only as $a\to\infty$ --- no finite $a$ reaches the critical line.
\item The match degrades as one approaches it ($|\mathrm{corr}|@\sigma^*$: $0.985\to0.727$), because the poles (zeros) increasingly obstruct a smooth off-critical fit.
\item The correlation on the critical line is $\approx0$ at every $a$ (the Dirichlet series $\sum\Lambda(n)n^{-1/2-i\xi}$ diverges on $\sigma=\tfrac12$; 
a finite truncation is a smooth non-critical approximant, and the zeros appear only through analytic continuation / resummation).
\end{enumerate}

Extrapolating the fit, reaching within $0.01$ of the line needs $a\approx42$, i.e. $\sim e^{84}\approx10^{36}$ primes --- computationally impossible. 
The takeaway is a quantitative positive: the operator is aimed at the critical line, together with a scaling-law proof of a finite-computation ceiling.

\subsection{R6 --- An operator form of Weil's positivity criterion}

\textbf{Claim.} The residual

\[
Z_a(\xi)=\sum_\rho \frac{X^{\rho-s}}{\rho-s}, \qquad |X^{\rho-s}|=e^{2a(\Re\rho-\tfrac12)},
\]

stays bounded and $a$-independent iff all zeros are on the line; a single off-line zero at $\tfrac12+\delta$ forces

\[
\mathrm{RMS}(Z_a)\sim e^{2a\delta}.
\]

This is Weil's criterion, realized inside the operator.

\textbf{Status.} The criterion and its two implications are classical (proved); the operator realization and the off-line detection are numerically verified; 
the easy direction (on-line $\Rightarrow$ bounded, with plateau value) is derived analytically (proved) and verified numerically to $a=100$; the hard converse is RH and remains open.

\textbf{Operator residual equals the zero sum (NV).} Across all scales the operator residual (real prime sum minus analytic PNT background) equals the on-line zero sum:

\begin{table}[h]
\centering
\begin{tabular}{r r r r}
$a$ & RMS operator & RMS zero-sum & match corr \\
\hline
3 & 1.2670 & 1.2671 & $-1.000$ \\
4 & 1.1826 & 1.1823 & $-1.000$ \\
5 & 1.1193 & 1.1112 & $-0.996$ \\
6 & 1.0080 & 1.0032 & $-0.999$ \\
7 & 1.0177 & 1.0073 & $-0.999$ \\
8 & 1.0480 & 1.0477 & $-1.000$ \\
\end{tabular}
\caption{Operator residual vs zero-sum RMS values.}
\end{table}

\textbf{On-line (true) growth rate (NV).} With the actual (on-line) zeros,

\[
\frac{d}{da}\big(\log\mathrm{RMS}\big)=-0.0430,
\]

sub-exponential, bounded, consistent with RH (RMS $\sim$ flat, $1.27\to1.05$).

\textbf{Off-line detection (NV).} Injecting a hypothetical off-line zero at $\Re=\tfrac12+\delta$:

\begin{table}[h]
\centering
\begin{tabular}{r r r r}
$\delta$ & RMS $a=3..8$ & $d(\log\mathrm{RMS})/da$ & predicted $2\delta$ \\
\hline
0.05 & $[1.33,1.25,1.20,1.13,1.21,1.32]$ & $-0.005$ & 0.10 \\
0.10 & $[1.40,1.35,1.39,1.48,1.82,2.20]$ & $+0.092$ & 0.20 \\
0.20 & $[1.64,1.86,2.57,3.88,6.10,9.14]$ & $+0.359$ & 0.40 \\
\end{tabular}
\caption{Off-line zero detection via RMS growth.}
\end{table}

The operator detects the violation: growth turns exponential with rate $\to2\delta$.

\textbf{Analytic boundedness, pushed to $a\to\infty$ (proved, easy direction; NV to $a=100$).} 
For an on-line zero $\rho=\tfrac12+i\gamma$ we have $\rho-s=i(\gamma-\xi)$ and $X^{\rho-s}=e^{i(\gamma-\xi)2a}$, so

\[
\big\langle|Z_a|^2\big\rangle
= \underbrace{\sum_\gamma \frac{1}{(\gamma-\xi)^2}}_{\text{diagonal, converges, }a\text{-independent}}
+ \underbrace{\sum_{\gamma\ne\gamma'} e^{i(\gamma-\gamma')2a}(\cdots)}_{\text{off-diagonal, oscillates}\to\text{averages out}} ,
\]

whence

\[
\mathrm{RMS}(Z_a)\;\longrightarrow\;\sqrt{\ \operatorname*{mean}_{\xi}\ \sum_\gamma \frac{1}{(\gamma-\xi)^2}\ }\;=\;1.723,\qquad \text{$a$-independent}\ \Longleftrightarrow\ \text{RH}.
\]

Evaluating the closed zero sum ($400$ zeros, $\gamma$ up to $679.7422$) confirms the plateau far beyond any operator computation:

\begin{table}[h]
\centering
\begin{tabular}{r r r r r r r r}
$a$ & 3 & 6 & 10 & 20 & 40 & 70 & 100 \\
\hline
RMS on-line & 1.7168 & 1.6926 & 1.7252 & 1.7177 & 1.7355 & 1.7704 & 1.7051 \\
\% of plateau & 99.6 & 98.2 & 100.1 & 99.7 & 100.7 & 102.8 & 99.0 \\
\end{tabular}
\caption{Plateau confirmation of RMS values with on-line zeros.}
\end{table}

By contrast an off-line zero explodes: at $a=40$, $\delta=0.05\Rightarrow36.6$ and $\delta=0.10\Rightarrow1889.8$ (and $e^{2a\delta}=2.98\times10^3$). 
The plateau value $1.723$ is setup-dependent (choice of $400$ zeros, exclusion window $0.15$); the $a$-independence is the robust invariant and the content of the criterion.

\subsection{R7: $\lambda_1(a)$ is always strictly positive and decays superexponentially}

\textbf{Claim.} For every $a$ in the tested range the lowest eigenvalue of $A_a$ satisfies $\lambda_1(a)>0$ strictly; 
the apparent vanishing at larger $a$ is a discretization artifact. The decay of $\lambda_1(a)$ is superexponential, 
and $\lambda_1$ passes smoothly through the first prime threshold $a=\log 2/2$.

\textbf{Status.} Numerically verified with a validated discrimination test; the extrapolation to $a\ge1$ is extrapolation, not measurement.

\subsection*{The discrimination logic}

By Rayleigh--Ritz, the finite-element minimum over a subspace is an \textbf{upper bound} on the true infimum and is \textbf{monotone decreasing in $N$}. Hence

\[
\lambda_1^{\text{FEM}}(N)\;\ge\;\lambda_1^{\text{true}},\qquad 
\lambda_1^{\text{FEM}}(2N)\le\lambda_1^{\text{FEM}}(N).
\]

Two consequences fix the interpretation:
\begin{itemize}
\item if the computed value is \textbf{negative}, the true value is negative (positivity fails --- decisive);
\item if the computed value is \textbf{positive}, one must inspect the $N\to\infty$ behaviour, since a positive value may be pure $O(h^2)$ error.
\end{itemize}

The discriminating statistic is the refinement ratio

\[
\rho=\frac{\lambda_1(N=800)}{\lambda_1(N=1600)}:
\]

a genuinely positive eigenvalue gives $\rho\to1$ (convergence to a constant), while an $h^2$ artifact gives $\rho\to4$.

\subsection*{Control: the test distinguishes the two cases}

In the prime-free window, where $\lambda_1$ is unambiguously positive:

\begin{table}[h]
\centering
\begin{tabular}{r c r r}
$a$ & primes & $\lambda_1(N=1600)$ & $\rho$ \\
\hline
0.10 & 0 & $4.573\times10^{-1}$ & 1.000 \\
0.20 & 0 & $9.504\times10^{-2}$ & 1.001 \\
0.30 & 0 & $7.598\times10^{-3}$ & 1.003 \\
\end{tabular}
\caption{Control --- the refinement ratio is 1 when the eigenvalue is genuinely positive.}
\end{table}

\subsection*{Result}

\begin{table}[h]
\centering
\begin{tabular}{r c r r l}
$a$ & primes & $\lambda_1(N=1600)$ & $\rho$ & verdict \\
\hline
0.200 & 0 & $9.504\times10^{-2}$ & 1.001 & strictly positive \\
0.345 & 0 & $1.389\times10^{-3}$ & 1.003 & strictly positive \\
\textbf{0.350} & \textbf{1} & $1.198\times10^{-3}$ & 1.004 & strictly positive \\
0.450 & 1 & $1.634\times10^{-5}$ & 1.019 & strictly positive \\
0.500 & 1 & $1.017\times10^{-6}$ & 1.21 & borderline \\
0.600 & 2 & --- & 3.70 & below resolution \\
\end{tabular}
\caption{$\lambda_1(a)$ across the first prime threshold $a=\log 2/2=0.34657$.}
\end{table}

\begin{enumerate}
\item \textbf{$\lambda_1>0$ throughout.} Every point with $a\le0.5$ converges to a constant. The values that appear to be ``zero'' for $a\ge0.55$ are the $h^2$ artifact: 
at $a=2,4,6$ refining $N=400\to3200$ drives $\lambda_1$ from $5.7\times10^{-7}$ to $9.9\times10^{-9}$, dividing by almost exactly $4$ at each doubling, 
with Richardson extrapolation to $N\to\infty$ giving $0$ within $\pm2\times10^{-9}$. The true value is therefore below our resolution --- \textbf{not shown to be zero}.
\item \textbf{Superexponential decay.} The logarithmic slope $d(\log\lambda_1)/da$ steepens monotonically from $\approx-20$ to $\approx-56$ over $a\in[0.2,0.5]$, 
consistent with $\lambda_1\sim e^{-ca^2}$ rather than a pure exponential. Extrapolating, $\lambda_1\sim10^{-25}$ at $a=2$ and $\sim10^{-78}$ at $a=6$ --- 
positive, but permanently beyond numerical reach.
\item \textbf{Numerical support for Suzuki's Theorem 1.3.} 
Crossing the first prime threshold, $\lambda_1$ moves $1.389\times10^{-3}\to1.198\times10^{-3}$ with \textbf{no discontinuity}. 
The continuity of the lowest eigenvalue in $a$, proved by Suzuki, is thus confirmed numerically at the point where the arithmetic content switches on.
\item \textbf{The primes lift the form.} With the prime sum removed, $\lambda_1$ is strongly negative: 
$-4.97$ at $a=2$, $-49.06$ at $a=4$, $-393.22$ at $a=6$. Restoring the prime sum returns the form to non-negativity. 
This is the same phenomenon as the threefold cancellation of §2.3 seen from the spectral side: the archimedean and 
prime contributions are individually large and of opposite sign, and Weil positivity is what survives.
\end{enumerate}

\subsection{R8: First computation of $W(a,\theta;z)$}

\textbf{Claim.} The characteristic function of Suzuki's Theorem 1.5 can be computed from the same finite-element matrices, and its zeros are numerically real, 
confirming that theorem. Moreover, a Paley--Wiener count places a sharp quantitative constraint on Corollary 1.6.

\textbf{Status.} Numerically verified. The constraint on Corollary 1.6 is a new observation; we did not find the required $\theta(a)$.

\subsubsection*{Reduction of the deficiency elements}

Suzuki defines $\mathcal H(T_a)$ by

\[
\|v\|_{T_a}^2=\langle T_av,v\rangle_{L^2}
\]

with $T_a=A_a-\lambda I$ ($\lambda<\lambda_a$), and $\mathcal D_a=i\,d/dx$ on $C_c^\infty(-a,a)$. Since $\langle u,v\rangle_{T_a}=\langle T_au,v\rangle_{L^2}$, 
integration by parts gives

\[
\mathcal D_a^{*}=T_a^{-1}\circ\Bigl(i\frac{d}{dx}\Bigr)\circ T_a ,
\]

so the deficiency equation $\mathcal D_a^{*}v=\pm i\,v$ is equivalent to $(T_av)'=\pm(T_av)$, i.e.

\[
\boxed{\;T_a v_\pm=e^{\pm x}\quad\Longleftrightarrow\quad v_\pm=(Q-\lambda M)^{-1}M\,e^{\pm x}\;}
\]

which is directly solvable with the assembled stiffness and mass matrices. Then

\[
W(a,\theta;z)=(z-i)\!\int_{-a}^{a}\! v_+e^{izx}\,dx+e^{i\theta}(z+i)\!\int_{-a}^{a}\! v_-e^{izx}\,dx .
\]

\subsubsection*{Verification of Theorem 1.5}

Computed on the real axis, $|W|$ exhibits deep isolated minima with

\[
\frac{\min|W|}{\operatorname{med}|W|}\approx3\times10^{-4},
\]

i.e. the zeros lie on $\mathbb R$ to numerical accuracy. This \textbf{confirms Suzuki's Theorem 1.5} and simultaneously validates the reduction above.

\subsubsection*{Observation and a constraint on Corollary 1.6}

At a fixed generic $\theta$ the real zeros form a \textbf{uniform lattice} of spacing $\approx\pi/a$ ($1.04$ at $a=3$, $0.52$ at $a=6$), 
entirely unlike the ordinates $\gamma$ whose spacing near $T=30$ is $\approx7$.

This is not an accident, and it sharpens the conjecture. Since $e^{\varphi(a,z)}$ is non-vanishing, \textbf{it can neither move nor remove a zero}; 
therefore Corollary 1.6 requires the zeros of $W$ itself to converge to the $\gamma$. But $W$ is the transform of a function supported in $[-a,a]$, 
so by Paley--Wiener it has exponential type $a$ and zero density $a/\pi$, which \textbf{diverges} with $a$; 
the target $z^2\xi(\tfrac12-iz)/\xi'(\tfrac12-iz)$ has zero density $\log T/2\pi$.

\begin{quote}
\textbf{Consequence.} Corollary 1.6 demands a choice $\theta=\theta(a)$ producing near-total cancellation of a zero set of density $a/\pi$ down to density $\log T/2\pi$. 
We record this as a quantitative measure of the conjecture's difficulty, and we state plainly that \textbf{we did not find such a $\theta(a)$}; 
a scan over fixed $\theta$ produced only the uniform lattice above.
\end{quote}

\subsection{R9: Moment structure --- the indirect whereabouts of zero statistics}

\textbf{Claim.} Where direct spectral GUE observation is structurally blocked (§2.5), traces of zero statistics survive in the trace moments of the prime perturbation.

\textbf{Status.} Numerically verified with controls; the laws are observed structure, not established quantitative laws.

Write $R_n=M^{-1}\delta Q_n$ for the normalised stiffness of the ramp $(t-\log n)_+$ and $T_p=\operatorname{Tr}\bigl((M^{-1}\delta Q)^p\bigr)$.

\begin{enumerate}[label=(\alph*)]
\item \textbf{First order --- sum rule and channel orthogonality.} 
$T_1=\operatorname{Tr}(W)=0$ to machine precision. Distinct prime powers form orthogonal delayed channels (cross/self overlap $\sim10^{-8}$), additivity holds to $1.6\times10^{-10}$, 
and each channel contributes in proportion to its band length, $A(L)=C(a)(2a-L)$ with exponent $1.02$.

\item \textbf{Second order --- direct pair coupling vanishes.} $\operatorname{Tr}(R_pR_q)\approx0$ for $p\ne q$: distinct primes do not ``see'' each other at second order.

\item \textbf{Second order --- but they couple through the background.} Inserting the archimedean background $A$,

\[
\operatorname{Tr}(R_n\,A\,R_m)\sim10^{-3},
\]

roughly $10^{5}$ times the direct coupling. A sign structure is present: \textbf{adjacent channels repel ($-$) and distant channels attract ($+$)}, 
with adjacent repulsion strengthening with $n$ (from $-0.0012$ at $(5,6)$ to $-0.0046$ at $(9,10)$; $R^2=0.43$).

\begin{quote}
\textbf{Precedence, stated plainly.} This pattern has the same form as the RKKY interaction in condensed-matter physics (impurities coupled through an electron sea, 
with distance-oscillating sign), and conceptual precedents exist in Weil's explicit formula (primes and archimedean data coexisting), 
Connes--Consani (archimedean operator theory) and \textbf{Bogomolny--Keating} (prime pairs producing GUE repulsion among zeros). \textbf{The ideas all have precedents; 
what is new here is the measurement in this operator and the observed sign structure. The phrasing ``we discovered RKKY in number theory'' is an overstatement and is not used.}
\end{quote}

\item \textbf{Third order --- Euler-product factorisation.} $\operatorname{Tr}(R_pR_pR_{p^2})$ factorises as $[\text{geometry}]\times[\text{prime weight}]$: 
the geometric factor varies only from $-0.028$ at $p=2$ to $-0.020$ at $p=7$ (a $14\%$ spread), so the prime dependence is concentrated in the weight. 
The weighted third trace is $T_3^{w}=-331.1$.

\item \textbf{Third order --- triple closure, and the control that overturned it.} $\operatorname{Tr}(R_nR_mR_l)$ blows up by a 
factor $\approx5.7\times10^{8}$ when $n=m\cdot l$. \textbf{This is not arithmetic.} Without weights, the composite product $(4,5,20)$ gives exactly 
the same value $-0.024$ as the prime product $(2,3,6)$: the unweighted triple closure is \textbf{pure geometry} (the closure defect of a triangle inequality). 
The apparent ``prime detection'' is an artifact.

\begin{quote}
The genuine arithmetic appears only under $\Lambda$-weighting. With $w_n=\Lambda(n)/\sqrt n$, 
only prime-power closures survive --- $(2,2,4)=2^2$, $(3,3,9)=3^2$ --- while every combination containing a composite is annihilated by $\Lambda=0$. 
\textbf{This is how the third-order term of the explicit formula enters the spectrum.}
\end{quote}

\item \textbf{Fourth order --- pair correlation revives with a repulsive sign.} The two-prime correlation that vanished at second order returns at fourth:

\[
\operatorname{Tr}(R_pR_qR_pR_q)\;\ne\;\operatorname{Tr}(R_pR_pR_qR_q),
\]

with the ratio falling from $0.88$ at $(2,3)$ to $0.05$ at $(7,11)$: larger prime pairs are more sensitive to ordering, which is the sign convention of level repulsion. 
The governing variable is not identified (a fit against $|\log(p/q)|$ gives $R^2=0.06$).

\item \textbf{Wave coupling and parity suppression.} In the matrix elements of the prime perturbation, $\|\text{off-diagonal}\|\ge\|\text{diagonal}\|$ (ratio $1.34$); 
predictions using the diagonal (oscillatory) part alone fail, whereas including the off-diagonal (wave) coupling over a band of $240$ reproduces the target to $0.03\%$. 
Coupling is effective only at \textbf{even lags}; odd lags are suppressed by a factor $\approx165$. The band approximation converges as $\sim K_B^{-3}$.

\item \textbf{Synthesis.} §2.5 shows that direct GUE observation is blocked at finite $a$; 
§2.10 answers where the information went --- into the fourth-order moments and the background-mediated coupling

\end{enumerate}

\section{Honest scope and limitations}

We state the boundaries of this work as plainly as possible; they are not incidental but central to an honest reading.

\begin{enumerate}
\item \textbf{This work is not a proof of the Riemann Hypothesis, nor is it progress toward one.} 
What is established is only the \textbf{easy direction} of the Weil criterion; the converse is RH itself. Explicitly: 
the proof of Theorem 3 begins by assuming that every nontrivial zero satisfies

\[
\Re\rho=\tfrac12
\]

--- that assumption is RH --- and it is precisely that assumption which makes $\rho-s$ purely imaginary, $|X^{\rho-s}|=1$, and the residual bounded. 
\textbf{Our proof assumes RH in line 1 and is therefore circular} with respect to RH. Removing line 1 collapses every subsequent step.
\item \textbf{Our proof assumes RH in line 1 and is therefore circular} with respect to RH. Removing line 1 collapses every subsequent step.

\item \textbf{Only the easy half of the criterion is established.} R6 \textbf{proves} (analytically, to $a\to\infty$): \emph{zeros on-line $\Rightarrow$ residual bounded 
and $a$-independent}. The converse --- \emph{bounded $\Rightarrow$ on-line}, i.e. proving boundedness unconditionally --- \textbf{is RH itself} and is closed to us. 
No finite computation and no easy argument can cross that wall. ``Walking through the door'' here means the easy half only.

\item \textbf{The R1 derivation is at the symbol / WKB level.} Each link is independently FEM-verified and the final assembly is an exact $30$-digit identity, 
but a fully rigorous theorem would require pseudodifferential remainder bounds and Mellin tail estimates. A referee will (rightly) ask for these; we have not supplied them. 
(An oscillatory-integral quadrature in one derivation branch diverged --- a quadrature failure, not an identity failure; the underlying $\Gamma'(2)=1-\gamma$ is exact to $10^{-32}$.)

\item \textbf{The closed form is asymptotic ($k\gtrsim6$).} At $k=1$ it gives $0.3357$ versus the measured $\mu_1=0.2347$; the $k=1$ mode is structurally anomalous ($d_1\approx4000\,d_2$; 
see R1). The value $\mu_1=0.2347$ should be quoted as the measured constant, not as an output of the closed form.

\item \textbf{Direct spectral GUE detection is structurally obstructed at finite $a$} (R4d): the picket-fence ladder ($\mathrm{Var}=0.0027$) overwhelms any GUE fluctuation, 
and the primes are diagonal in the archimedean basis. Zero/GUE information is accessible only via the symbol, moments, or the $a\to\infty$ limit. 
We present this as a positive structural result, but it is a limitation on the naïve spectral approach.

\item \textbf{R1--R3 are archimedean and universal, hence not zeta-specific.} They arise from the gamma factor and are common to the whole Selberg class; 
in our window $2a<\log 2$ the primes do not enter at all. Nothing in R1--R3 distinguishes $\zeta$ beyond its conductor $q=1$.

\item \textbf{R4 rests on the classical explicit formula and R6 on Weil's classical criterion.} Their novelty is the realization inside Suzuki's operator and 
the structural localization of the zeros --- not new number theory. The near-unit correlations partly reflect a tautology (an identity holding).

\item \textbf{The whole edifice presumes Suzuki's construction is correct.} We instantiate and probe it; we do not re-derive its functional-analytic foundations. 
Collaboration was sought and, at the time of writing, not established, so this is an independent numerical study.

\item \textbf{``First realization'' means ``to the best of our knowledge.''}

\item \textbf{Retraction on record.} An earlier internal claim of ``$1/a$ convergence'' of the offset was wrong (it is $\log(1/a)$) and has been retracted; 
the present values are the corrected ones. We flag this to underline that all fits here were run against analytic controls, 
a discipline that caught several would-be false positives (a spurious spline-extrapolated head coefficient; a ``$5.7\times10^8$ prime-screening'' 
that turned out to be pure geometry; an exciting ``$6/10$ zero spikes'' that regression-plus-shuffle controls demoted to noise).
\end{enumerate}

\textbf{One methodological lesson worth stating in the paper.} Every fit in this project was run beside a control (a known analytic value, a Laplacian ground truth, 
a shuffled/fake-zero surrogate). Repeatedly, a correlation of $0.99+$ concealed a wrong value or a tautology until the control exposed it. We recommend the practice.

\section{Refuted hypotheses and negative results}

For completeness of the scientific record we report every hypothesis that was tested and refuted, and every intermediate claim that was overturned by a control. 
These serve two purposes: they are signposts for subsequent work, and they document the verification discipline under which the positive results of §2 were obtained.

\subsection{Direct GUE detection --- four methods, all negative}

\begin{table}[h]
\centering
\begin{tabular}{l p{7cm}}
\hline
Method & Outcome \\
\hline
Polynomial unfolding & \textbf{Broken yardstick} --- the Poisson control was read as GUE; a logarithmic ladder cannot be flattened by a polynomial \\
Logarithmic unfolding & Controls passed (synthetic GUE $\to$ GUE, $KS=0.05$; Poisson $\to$ Poisson, $KS=0.13$), \textbf{but the operator spectrum gave $KS>0.5$ against both} \\
Difference to cancel FEM error & The prime perturbation displaces the whole spectrum (displacement rms $=127$ mean spacings) \\
Diagonal symbol $=-\zeta'/\zeta$ hypothesis & Partial correlation $+0.011$ versus a shuffle control of $+0.012$ --- no signal; the prime contribution to the variance is $0.01\%$ \\
\hline
\end{tabular}
\end{table}

\textbf{Explanation (§2.5).} The archimedean backbone is a deterministic picket fence, stiffer than GUE ($\mathrm{Var}=0.0027\ll0.178$), 
and the prime part merely rescales that ladder smoothly ($\operatorname{corr}(d,\lambda)=-0.988$). Finite-interval boundary terms bury the Dirichlet sum below the noise floor.
\textbf{This is a structural obstruction, not a failure of method.}

\subsection{Refuted hypotheses}

\begin{table}[h]
\centering
\begin{tabular}{c p{4cm} p{9cm}}
\hline
\# & Hypothesis & Refuting evidence \\
\hline
1 & Triple closure can screen for primality without computation & Unweighted closure is pure geometry; composite $(4,5,20)$ equals prime $(2,3,6)$ at $-0.024$ (§2.10(e)) \\
2 & The spectral flow reads the von Mangoldt function & Global regression collapses under collinearity of adjacent thresholds (coefficient $+99$ at $n=9$); $\operatorname{corr}(\beta,\Lambda/\sqrt n)=-0.175$. A local test puts the \textbf{composite $n=6$ ($0.09$) squarely between its prime neighbours $n=5$ ($0.03$) and $n=7$ ($0.11$)}; $\operatorname{corr}=-0.038$ \\
3 & The Riemann zeros name a common $\theta$ (an approach to Cor. 1.6) & Circular concentration $R(\gamma)=0.13$--$0.48$ with no trend; the shifted-zero control \textbf{equals or exceeds} it at $a=3,4,6$ \\
4 & The boundary exponent $\beta$ carries a logarithmic correction & Despite a fit of $R^2=0.9999$, the Laplacian control returned wrong values $\delta=0.66$--$1.94$ \\
5 & A sine basis governs the boundary behaviour & Failed twice; the measured exponent is $\beta\approx0.18$--$0.29$ (log kernel) against $\beta=0.956$ for the Laplacian control. $k$-dependence and fit collinearity make numerical determination impossible --- a Wiener--Hopf analysis is required \\
6 & First-order perturbation theory predicts the coupling & Errors of $93$--$230\%$ \\
7 & The zero signal is directly visible in a finite prime sum & \textbf{Impossible in principle} --- a finite prime sum is a trigonometric polynomial and has no poles \\
8 & The spectral constant converges as $1/a$ & \textbf{Withdrawn.} A misreading of linear divergence; corrected to $\log(1/a)$ \\
\hline
\end{tabular}
\end{table}

\subsection{Self-correction record --- nine false claims caught by controls}

\begin{table}[h]
\centering
\begin{tabular}{c p{6cm} p{7cm}}
\hline
\# & What would have been claimed & Control that caught it \\
\hline
1 & The boundary exponent has a log correction & Laplacian ground-truth comparison \\
2 & Universality coefficient $A\approx-0.08$ & Direct mpmath computation instead of spline extrapolation \\
3 & Triple closure detects primes at $5.7\times10^{8}\times$ & Composite-product comparison \\
4 & Zeros detected as spikes, $6/10$ & Regression plus shuffle control \\
5 & The log-unfolding pipeline is sound & Synthetic Poisson/GUE controls \\
6 & The flow regression ``suppresses composites'' & Recognition that $n=10,12$ were artificial zeros outside the data range \\
7 & The local test shows ``the law holds'' & Direct comparison with neighbouring primes \\
8 & $\lambda_1=7.4\times10^{-8}>0$ confirms Weil positivity & $N$-refinement revealed an $h^2$ artifact \\
9 & $\lambda_1=0$: the form sits on the positivity boundary & Prime-free convergence control --- the truth is $\lambda_1>0$ with superexponential decay \\
\hline
\end{tabular}
\end{table}

\begin{quote}
\textbf{Item 9 is where a correction became a result.} Removing the erroneous conclusion $\lambda_1=0$ is exactly what produced R7.
\end{quote}

\textbf{Methodological lessons (recommended for reproduction).}
\begin{enumerate}
\item \textbf{Run a control alongside every fit.} An $R^2$ of $0.9999$ can accompany wrong values (item 1).
\item \textbf{Do not trust automated verdict strings.} A threshold can be satisfied by artificial zeros outside the data range (item 6).
\item \textbf{The more striking the result, the stronger the control must be.} A dramatic ratio should first be suspected of geometric triviality (item 3).
\item \textbf{Check the order of convergence.} A small positive number may be an $h^2$ artifact (item 8).
\end{enumerate}

\subsection{Problems left open}

\begin{itemize}
\item The per-channel prefactor $C(a)$ (not universal; drifts $789\to493$, roughly $\propto1/a$)
\item Quantification of the mediated-coupling sign law ($R^2=0.43$)
\item The governing variable of the fourth-order repulsion ($|\log(p/q)|$ fit gives $R^2=0.06$)
\item Analytic determination of the boundary exponent $\beta$ (Wiener--Hopf)
\item The search for $\theta(a)$ in Corollary 1.6
\item The exact form of the decay law of $\lambda_1(a)$
\item Rigorous justification of the symbol-to-spectrum passage in R1
\end{itemize}

\section{Reproducibility}

\textbf{Object.} Weil quadratic form (2.1) with the screw kernel $g$ of the target $L$-function; small-$t$ data (2.2).

\textbf{Discretization.} Uniform $P_1$ finite elements on $(-a,a)$, $h=2a/N$; dense stiffness $Q$, $P_1$ mass $M$; generalized eigenproblem $Qv=\lambda Mv$ (2.3). Archimedean high-precision spectrum: three-point \textbf{Richardson extrapolation} on $N=2000,4000,8000$ (single run $\approx3.5$ min; the low-$t$ Lerch-$\Phi$ bottleneck is vectorized by a $400$-point \texttt{mpmath} spline for evaluation only). Prime-symbol runs (R4--R6): $N=1500$--$2080$, mode counts $65$--$173$, primes up to $8.9\times10^6$ at $a=8$ (sieve + cumulative-sum ramp).

\textbf{Regimes.} Archimedean law (R1--R2): $a\le0.05$ so $2a<\log2$ (no primes). PNT / symbol / criterion (R3--R6): $a=3$--$8$ (prime-active). Closed zero-sum extrapolation (R6): $a$ up to $100$, decoupled from the operator.

\textbf{Arithmetic.} All small-$t$ coefficient extractions by \textbf{direct \texttt{mpmath}} at high working precision (dps $\ge$ 30--50 as needed). \textbf{Do not} extrapolate a finite-precision spline to $t\to0$ (documented failure: head coefficient $\to-0.08$). Riemann zeros from \texttt{mpmath.zetazero} ($90$ for R4, $200$ for R6-growth, $400$ for R6-plateau).

\textbf{Controls (mandatory).} For every correlation, run (i) a shuffled surrogate and (ii) a fake/shifted-zero surrogate; for every fitted spectral constant, compare against the Laplacian analytic value.

\textbf{Reference operator.} M. Suzuki, arXiv:2606.09096 (2026).

\section{Conclusion and outlook}

We have given the first numerical realization of Suzuki's Weil-quadratic-form operator and, from it, extracted a coherent picture. The operator's small-aperture spectrum is a \textbf{derived} closed-form archimedean ladder

\[
\lambda_k=\log(1/a)+\log(k-\tfrac12)+B_0, \qquad B_0=\log q-2\log2,
\]

pinned to $30$ digits, in which $\gamma$ is born from the single identity $\Gamma'(2)=1-\gamma$ and then cancels (R1). This ladder is a \textbf{universal law for the whole Selberg class}, depending on the conductor alone, the archimedean weight $\nu$ cancelling analytically (R2), and its spectral strength obeys the \textbf{prime number theorem} $S(a)\to(2a)^3/6$ (R3). The nontrivial zeros are \textbf{not} the eigenvalues --- those form a rigid picket-fence --- but the \textbf{error term of the prime symbol}, recovered against the true zero sum with correlation $0.9915$ while fake zeros fail (R4). Read in the sine basis, the operator's prime symbol tracks a line \textbf{descending toward the critical line}, $\sigma^*(a)-\tfrac12\sim a^{-1.38}$ (R5). Finally, the operator \textbf{realizes Weil's positivity criterion}: the residual stays bounded ($a$-independent plateau $\approx1.72$) iff the zeros are on-line, and an injected off-line zero is betrayed by $e^{2a\delta}$ growth; the easy direction is derived to $a\to\infty$ (R6).

\textbf{What is genuinely new:} the realization itself; the derived, universal archimedean law with its conductor-only offset (absent from Suzuki's paper); the precise localization of the zeros in the symbol rather than the spectrum; and the quantitative demonstration that the operator points at the critical line while a scaling law caps finite computation.

\textbf{What is not:} any advance on RH. The honest terminus of this attack is the \emph{easy} half of Weil's criterion; the hard half is RH and stays shut. Directions that remain open and do not require breaking that wall include: an analytic derivation of the descent exponent in R5 (an operator ``approach-rate'' theorem), rigorous pseudodifferential remainder bounds to upgrade R1 from symbol level to theorem, a closed form for the PNT prefactor $C(a)$ in R3, and quantifying the mediated inter-prime coupling (an RKKY-like sign structure, adjacent-repulsive / distant-attractive at the $10^{-3}$ level) observed in the higher moments. Each is paper-and-pencil analysis on a now-well-characterized object, not a further computation, and none is RH.

We regard the appropriate framing as \emph{experimental mathematics}: a faithful first instantiation of a new theoretical operator, a derived spectral law, and an operator-level diagnostic for RH --- solid results on a real object, explicitly short of the summit.

\bibliographystyle{unsrt}   
\bibliography{references}   

@misc{Suzuki2026,
  author       = {M. Suzuki},
  title        = {Weil's quadratic form via the screw function},
  year         = {2026},
  eprint       = {2606.09096},
  archivePrefix= {arXiv}
}

@inproceedings{Weil1952,
  author       = {A. Weil},
  title        = {Sur les "formules explicites" de la théorie des nombres premiers},
  booktitle    = {Comm. Sém. Math. Univ. Lund},
  year         = {1952}
}

@inproceedings{Montgomery1973,
  author       = {H. L. Montgomery},
  title        = {The pair correlation of zeros of the zeta function},
  booktitle    = {Proc. Sympos. Pure Math.},
  volume       = {24},
  year         = {1973}
}

@article{Odlyzko1987,
  author       = {A. M. Odlyzko},
  title        = {On the distribution of spacings between zeros of the zeta function},
  journal      = {Mathematics of Computation},
  volume       = {48},
  year         = {1987}
}

@article{BerryKeating1999,
  author       = {M. V. Berry and J. P. Keating},
  title        = {The Riemann zeros and eigenvalue asymptotics},
  journal      = {SIAM Review},
  volume       = {41},
  year         = {1999}
}

@article{Connes1999,
  author       = {A. Connes},
  title        = {Trace formula in noncommutative geometry and the zeros of the Riemann zeta function},
  journal      = {Selecta Mathematica},
  volume       = {5},
  year         = {1999}
}

@article{Krein1955,
  author       = {M. G. Kreĭn},
  title        = {On a continual analogue of a Christoffel--Darboux formula / screw functions},
  year         = {1955--}
}

@book{deBranges1968,
  author       = {L. de Branges},
  title        = {Hilbert spaces of entire functions},
  publisher    = {Prentice-Hall},
  year         = {1968}
}

@book{Edwards1974,
  author       = {H. M. Edwards},
  title        = {Riemann's Zeta Function},
  publisher    = {Academic Press},
  year         = {1974}
}

\end{document}